\documentclass[a4paper]{article}
\pdfoutput=1
\usepackage{hyperref}
\hypersetup{
	pdfinfo={
		Title={Learning based Hamilton–Jacobi–Bellman Methods for Optimal Control},
		Author={Sixiong You, Ran Dai, and Ping Lu},
		keyward={Supervised Learning; Reinforcement Learning; Optimal Control; Two Point Boundary Value Problem, Hamilton–Jacobi–Bellman Equation}
	}
}
\usepackage{pdfpages}
\usepackage{amsmath}
\usepackage{amsfonts}
\usepackage{hyperref}
\usepackage{graphicx}
\usepackage{subcaption}
\usepackage{blindtext}
\usepackage{epsfig} 

\begin{document}

\center{ \bf
	Learning based Hamilton–Jacobi–Bellman Methods \\for Optimal Control}

\center{Sixiong You, Ran Dai, and Ping Lu}


\begin{abstract}

Many optimal control problems are formulated as two point boundary value problems (TPBVPs) with conditions of optimality derived from the Hamilton–Jacobi–Bellman (HJB) equations. In most cases, it is challenging to solve HJBs due to the difficulty of guessing the adjoint variables. This paper proposes two learning-based approaches to find the initial guess of adjoint variables in real-time, which can be applied to solve general TPBVPs. For cases with database of solutions and corresponding adjoint variables of a TPBVP under varying boundary conditions, a supervised learning method is applied to learn the HJB solutions off-line. After obtaining a trained neural network from supervised learning, we are able to find proper initial adjoint variables for given boundary conditions in real-time. However, when validated solutions of TPBVPs are not available, the reinforcement learning method is applied to solve HJB by constructing a neural network, defining a reward function, and setting appropriate super parameters. The reinforcement learning based HJB method can learn how to find accurate adjoint variables via an updating neural network. Finally, both learning approaches are implemented in classical optimal control problems to verify the effectiveness of the learning based HJB methods.

\end{abstract}

\section{INTRODUCTION}

Optimal control is to find a control policy for a dynamic system such that the cost function is minimized subject to certain constraints \cite{bryson2018applied}. Due to the complexity of real-world dynamic systems, 
many optimal control problems involve nonlinear dynamics and/or nonconvex constraints. Existing literatures for solving optimal control problems are generally classified into two categories, direct methods and indirect methods \cite{betts1998survey}.
%
%

The basic idea of direct methods is to approximate an optimal control problem via well-defined collocation methods, e.g., \cite{fahroo2002direct}, to convert it into a parameter optimization problem which can be solved by the state-of-the-art nonlinear optimization methods, such as nonlinear programming \cite{bertsekas1997nonlinear}. 
%
%
%
%
Other the other hand, indirect methods find the sufficient and necessary optimality conditions based on the Hamilton–Jacobi–Bellman (HJB) partial differential equations \cite{bryson2018applied}. However, only a few special cases can find the analytical solution of HJB equations due to nonlinearity. 
There are two approaches to solve nonlinear HJBs. One is to utilize numerical methods to find the approximate solutions of HJBs, such as Galerkin approximation \cite{Randal1997Galerkin}, recursive approximation \cite{George1979Approximation, Beard1998Approximate}, asymptotic approximation \cite{Evans2017Approximating}, just to name few. The other one is to search the adjoint variables to integrate the partial differential equations, such as the single shooting and multiple shooting methods \cite{bock2000direct}. However, due to the sensitivity of adjoint variables, existing methods cannot guarantee convergence to the exact adjoint variables, which makes the HJB based method infeasible for general optimal control problems.   

This paper revisits the HJB based method by effectively searching the initial adjoint variables to solve general optimal control problems in real-time. To achieve this goal, machine learning is adopted to improve the effectiveness and efficiency of searching initial adjoint variables. Machine learning is a branch of artificial intelligence and has been successfully applied in many fields \cite{Kotsiantis2006Machine, ALJARRAH2015BigData, Bengio2013Representation}. In general, machine learning can be classified as supervised learning, unsupervised learning, and reinforcement learning. For supervised learning, all training data is labeled, which means all inputs have corresponding outputs, and then it will train the neural network to achieve the least total error between predictions and labels \cite{Kotsiantis2007Supervised}. If the objective is to classify the inputs into different groups based on the similarity among them and there is no label for the training data, then the unsupervised learning is able to solve this type of problems \cite{Yoshua2012Unsupervised}. Instead of learning from the training database, the reinforcement learning is to learn via interactions with the environment and simultaneously optimize the policy by iteratively choosing actions \cite{kiumarsi2018optimal}. 
Reinforcement learning has been applied to search for an approximate HJB solution \cite{bohme2017hybrid, Tamimi2008Discrete, Kiumarsi2014Reinforcement, LUO2015Reinforcement}. Among these works, most are designed for the discrete-time systems. For continuous control problems, it requires further investigation to reduce sample complexity, effectively use available data, and properly choose a network to approximate the value function.

Considering the advantages of solving HJBs using a good guess of adjoint variables, especially for continuous control problems with high precision requirement on solutions, we focus on searching for the initial adjoint variables. Two learning based HJB approaches are proposed, including supervised learning based Hamiltonian (SLH) and reinforcement learning based Hamiltonian (RLH). Since there are many numerical optimization approaches that can solve HJB equations off-line via approximate formulation and/or different initial guess, it is practical to construct the training database with correct solutions of a specific optimal control problem under varying boundary conditions. The supervised learning is able to learn the HJB solutions off-line and the trained neural network will find the proper adjoint variables in real-time for the corresponding problem with given boundary conditions. The combined supervised learning and HJB lead to the SLH method. However, for challenging optimal control problems that cannot be easily solved via off-line optimization methods to construct the training database, the RLH method is proposed aiming to learn the adjoint variables from experience in the absence of reliable training database.


This paper is organized as follows. In Section II, the general formulation of optimal control problem and the HJB equations are introduced. In Section III, the proposed learning based HJB approaches are presented. In Section IV, the results of these two learning based HJB methods for solving optimal control problems are shown and analyzed. Conclusions and future work are presented in Section V.
\section{OPTIMAL CONTROL AND HJB EQUATIONS}

The optimal control problem of a continuous-time system can be expressed as
\begin{eqnarray}\label{eq:OCP1}
&\min\limits_{\mathbf{u}} & J = \int_{0}^{t_f} g(\mathbf{x},\mathbf{u},t)dt \\
&\text{s.t.}  & \dot{\mathbf{x}} = f(\mathbf{x},\mathbf{u},t) \nonumber\\
&             & \mathbf{x}(0) = \mathbf{x}_0,\, \mathbf{x}(f) = \mathbf{x}_f \nonumber
\end{eqnarray}
where $\mathbf{u}$ is the control variable, $\mathbf{x}$ is the state vector, $t$ represents time, $J$ is the objective function, $\dot{\mathbf{x}} = f(\mathbf{x},\mathbf{u},t)$ represents the system dynamics, $\mathbf{x}_0$ and $\mathbf{x}_f$ represent initial and final boundary values, respectively, and $t_f$ refers to the terminal time. 
By introducing the adjoint variable set $\lambda$, the Hamiltonian is defined as
\begin{eqnarray}\label{eq:H1}
H(\mathbf{x},\mathbf{u},t) = g(\mathbf{x},\mathbf{u},t)+ \lambda  f(\mathbf{x},\mathbf{u},t). 
\end{eqnarray}
By deriving the first-order optimality conditions of (\ref{eq:H1}), the HJB equations in form of partial differentials are written as
\begin{eqnarray}\label{eq:H2}
\dot{\mathbf{x}}^* = \frac{\partial \mathbf{H^*}}{\partial \lambda^*},\, \dot{{\lambda}}^* = \frac{\partial \mathbf{H^*}}{\partial \mathbf{x^*}}
\end{eqnarray}
where $\mathbf{^*}$ represents corresponding values at optimum point with boundary conditions $\mathbf{x}(0) = \mathbf{x}_0,\, \mathbf{x}(f) = \mathbf{x}_f$.

If the values of adjoint variables, $\lambda^*$, at the initial point $t=0$ are given, we are able to find the optimal solution in continuous time by integrating $\dot{\mathbf{x}}^*$ and $\dot{{\lambda}}^*$ forward to any desired time interval. However, the adjoint variables have no physical meaning in an optimal control problem, it is challenging to guess their values, even their value range. In addition, for nonlinear differentials, $\frac{\partial \mathbf{H^*}}{\partial \lambda^*}$ and $\dot{{\lambda}}^* = \frac{\partial \mathbf{H^*}}{\partial \mathbf{x^*}}$, it is not feasible to find their analytical solutions in general cases. Although multiple shooting methods \cite{bock2000direct} and nonlinear optimization algorithms \cite{bertsekas1997nonlinear} have been developed to solve HJBs, they cannot guarantee yielding an optimal solution in real-time for general cases, especially for systems with highly nonlinear dynamics.   
Therefore, we propose the learning based HJB methods focusing on finding proper adjoint variables at the initial point of problem \eqref{eq:OCP1}, which leads to a continuous optimal control law in real-time.

\section{LEARNING BASED HJB METHODS}

In this section, two learning based HJB methods, including SLH and RLH, are introduced to solve optimal control problems in real-time. Firstly, when the off-line solutions of a specific problem under varying boundary conditions are available to construct the training database, a supervised learning algorithm combined with HJB is developed, which can be applied to general optimal control problems, even for cases with hypersensitive adjoint variables. Next, a deep reinforcement learning based HJB approach is introduced which can efficiently solve the optimal control problems without validated training data. 

\subsection{Supervised Learning based Hamiltonian (SLH)}
Taking the boundary conditions and terminal time as an input set $i$, denoted as $\mathbf{x}^i_{\text{input}} = [\mathbf{x}_0^i,\mathbf{x}_f^i,t_f^i]$, and the corresponding adjoint variables at initial time as the output, denoted as $\mathbf{x}^i_{\text{label}} = \lambda_0^i$, with sufficient database recording inputs and outputs, the initial adjoint variables can be found from the mapping function between $\mathbf{x}_{\text{input}}$ and $\mathbf{x}_{\text{label}}$. To obtain the mapping function, we employ the artificial neural network (ANN) to map the relationship between input and output. 

\begin{figure}[!h]
	\centering
	\includegraphics[width=8.6cm]{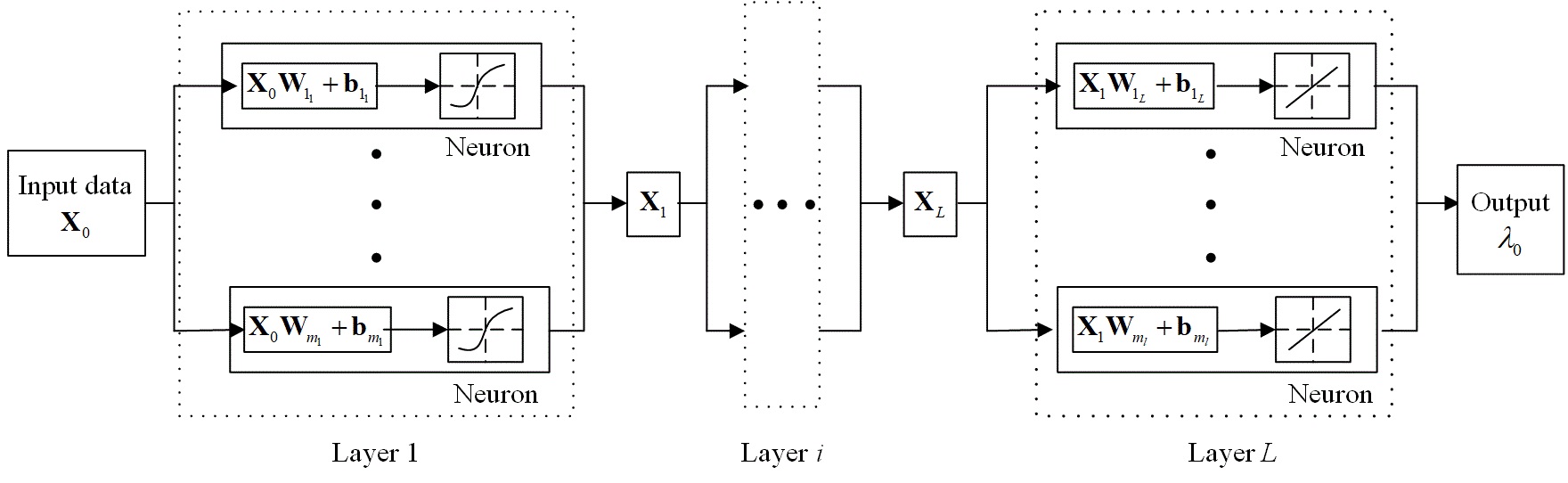} 
	\caption{The network architecture of SLH}
	\label{f:ANN}
\end{figure}

The constructed network architecture is shown in Fig. \ref{f:ANN}, where $L$ fully connected layers are included and each layer has $m_l$ neurons except the output layer. The number of neurons in the output layer is the same as the number of elements in $\lambda_0$. To better fit the nonlinearity of different systems, the activation functions $\mathbf{F_{act}}$ can be adjusted independently. By constructing a suitable network, the SLH is developed to find $\lambda_0$. Combining the given initial states, $\mathbf{x}_0$, HJB equations in (\ref{eq:H2}) can be integrated forward to any desired time interval to obtain the optimal solutions in real-time.

However, there are a set of optimal control problems with very sensitive adjoint variables. Very small changes, e.g., less than $10^{-4}$ magnitude, of the adjoint variables may lead to significant difference of the final solution. These optimal control problems are named hypersensitive HJBs where the dynamics exhibit fast contraction and expansion along the time interval \cite{aykutlug2016manifold}. Consequently, the accuracy requirement of the initial adjoint variables $\lambda_0$ may go beyond the precision scale provided from a constructed neural network. Thus, searching for $\lambda_0$ only without guarantee of high precision cannot always lead to an optimal solution for the hypersensitive cases. Generally, the hypersensitive HJBs can be further divided into completely hypersensitive HJBs and partially hypersensetive HJBs, and here we will focus on the completely hypersensitive HJBs. 
%

It has been found that the solution of a completely hypersensitive HJB can be approximated by three phases along the time, including stable phase, unstable phase, and equilibrium point \cite{aykutlug2009approximate}, expressed as 
\begin{eqnarray}\label{eq:S1}
\mathbf{p}=\begin{cases}
\mathbf{p}_s(t), & 0\leq t \leq t_{ib}\\
\mathbf{p}_e, & t_{ib} < t < t_{fb}\\
\mathbf{p}_u(t), & t_{fb} \leq t \leq t_{f}
\end{cases}
\end{eqnarray} 
where $\mathbf{p}=[\mathbf{x}^T,\lambda^T]^T$, $\mathbf{p}_s(t)$ is the value of $\mathbf{p}$ in the stable phase, and $\mathbf{p}_u(t)$ is the value of $\mathbf{p}$ in the unstable phase, and  $\mathbf{p}_e$ is the equilibrium point of $\mathbf{p}$ in the phase space, $t_{ib}$ and $t_{fb}$ are the ending time of stable phase and starting time of unstable phase, respectively. 

\begin{figure}[!h]
	\centering
	\includegraphics[width=8.8cm]{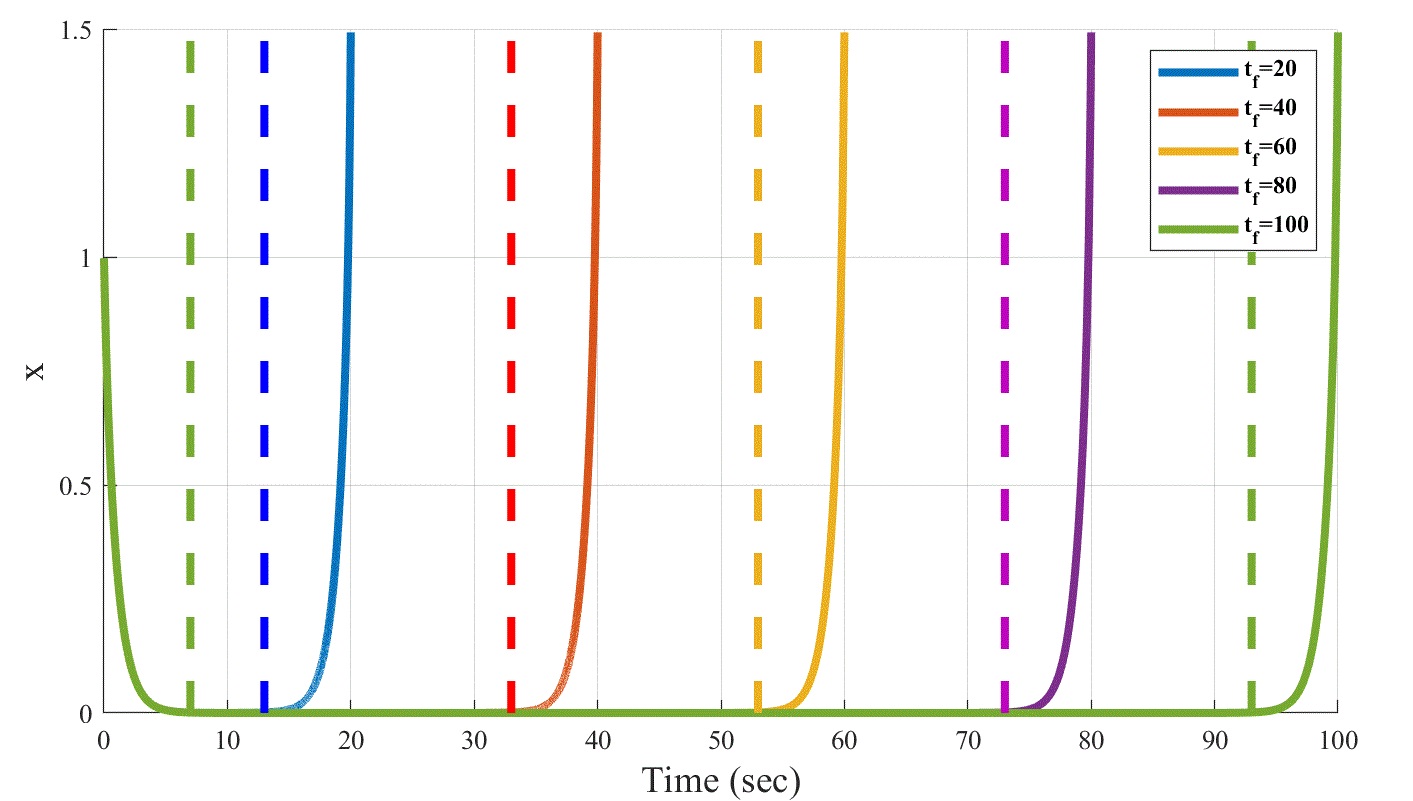} 
	\caption{An examples of completely hypersensitive HJBs}
	\label{f:HHJB}
\end{figure}

To illustrate it more clearly, an example of completely hypersensitive HJBs is presented here, which is expressed by
\begin{eqnarray}\label{eq:M1}
&\min\limits_{\mathbf{u}} & \mathbf{J}_2 = \int_{0}^{t_f} (x^2+u^2)dt \\
&\text{s.t.}  & \dot{{x}} = -x^3 + u \nonumber\\
&            & \mathbf{x}(0)=1,\, \mathbf{x}(t_f)=1.5 \nonumber
\end{eqnarray}
The optimal solution of this hypersensitive example with varying terminal time is shown in Fig. \ref{f:HHJB}. When the terminal time is sufficiently large, e.g., $t_f>20$ sec, the optimal solution can be approximated by three phases, named stable phase, equilibrium point, and unstable phase,  with phase separation marked by dashed lines. 

To solve the completely hypersensitive HJBs, an extended SLH method is proposed to learn multi-phase adjoint variables. At the equilibrium point in \eqref{eq:S1}, $\mathbf{p}_e$ is a constant number. So only the stable phase and unstable phase have time varying states and adjoint variables. For the stable and unstable phases, each phase can be further divided it into several subsegments with associated ``initial" adjoint variables for each subsegment, aiming to lower the sensitivity. To develop a systematic approach to solve the hypersensitive HJBs using the SLH framework, the first step is to determine the three phases based on available HJB solutions with varying boundary conditions. After selecting the transition points between stable phase and unstable phase, each phase is further divided into several subsegments. The SLH method can then be applied to learn the ``initial" adjoint variables of every subsegment off-line. Finally, the hypersensitive HJB can be solved in real-time by interpolating the trained neural network on-line for given boundary conditions.

\subsection{Reinforcement Learning based Hamiltonian (RLH)}

For challenging optimal control problems with highly nonlinear dynamics, existing approaches may not able to find their optimal solutions off-line, which means there are no enough effective data to train a neural network. To solve these problems, the deep reinforcement learning will be applied, which will auto tune the neural network to find the optimized adjoint variables.

\begin{figure}[!h]
	\vspace{.2cm}
	\centering
	\includegraphics[width=6.5cm]{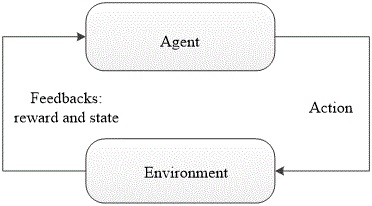} 
	\caption{The basic composition of reinforcement learning}
	\label{f:RL}
	\vspace{-.2cm}
\end{figure}

The basic framework of reinforcement learning is shown in Fig. \ref{f:RL}, which includes two major components, agent and environment. The agent will determine the policy of taking actions $\mathbf{a_t}$ for different input states $\mathbf{s}_t$, where $\mathbf{*}_t$ means $\mathbf{*}$ at stage $t$. The environment will give feedbacks for different input actions $\mathbf{a_t}$. In General, the feedbacks include the corresponding rewards $\mathbf{R}_t$ and next states $\mathbf{s}_{t+1}$. For the studied HJB, we define state, action, and rewards as
\begin{eqnarray}\label{eq:R1}
& \mathbf{a_t} = \lambda_0 \\ \nonumber
& \mathbf{s_t} = [\mathbf{x_0},\mathbf{x_f},t_f,...] \\\nonumber
& \mathbf{R_t} = h(\mathbf{a_t},\mathbf{s_t}) \nonumber
\end{eqnarray}
where $\lambda_0$ is the adjoint variables at initial time. In addition, as only one step is required to solve HJB, the next stage will be the terminal stage, the corresponding $\mathbf{R_t}$ can be accurately obtained without estimating $\mathbf{R_{t+1}}$, which makes the learning process more effective. 

According to the type of action space, the existing reinforcement learning algorithms can be divided into two categories. For the discrete action space, examples of developed methods include Monte Carlo \cite{Evans1998Introduction}, Q-learning \cite{Watkins1992Q}, State-action-reward-state-action \cite{Rummery1994On-line}, and deep Q network \cite{Volodymyr2015Human}. For continuous action space, examples include asynchronous advantage actor-critic algorithm \cite{Volodymyr2016Asynchronous}, deep deterministic policy gradient (DDPG) \cite{Timothy2015Continuous} and proximal policy optimization \cite{Schulman2017Proximal}. Since the HJB problem has continuous action space and deterministic environment, the DDPG algorithm is adopted here.    
\subsubsection{Introduction of DDPG Algorithm}
\begin{figure}[!h]
	\vspace{-.2cm}
	\centering
	\includegraphics[width=8.5cm]{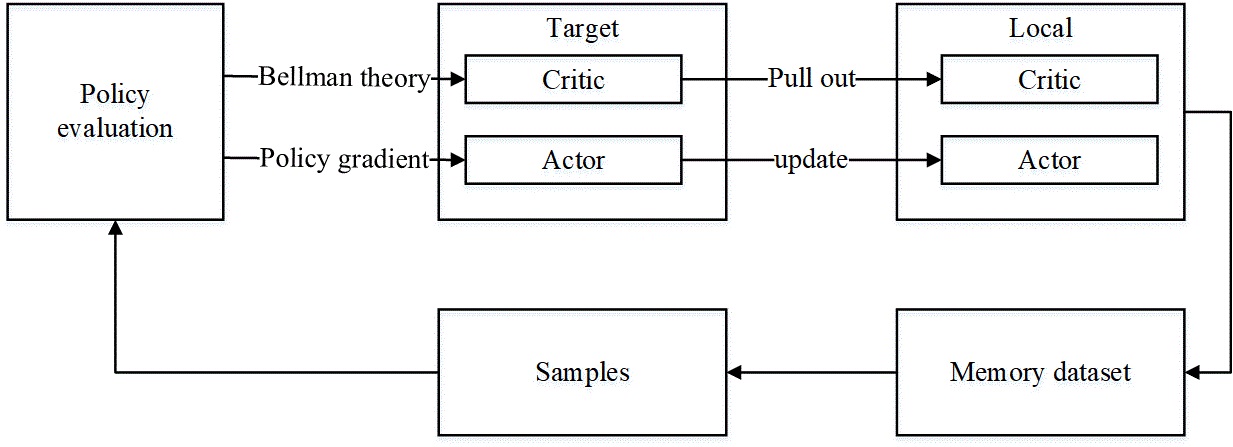} 
	\caption{The framework of DDPG}
	\label{f:DDPG}
	\vspace{-.3cm}
\end{figure}

DDPG can be regarded as the combination of actor-critic algorithm and deep Q learning, where the neural networks is divided into two levels, local network and target network. In each episode, the target network is softly updated, which means it will be updated slowly, and the local network will copy the target network after fixed episodes. In each level, there are two neural networks, actor network and critic network. The actor network can output actions $\mathbf{a_t}$ according to the observation states $\mathbf{s_t}$, and the critic network can output the reward value $\mathbf{R_t}$ based on the input states $\mathbf{s_t}$ and actions $\mathbf{a_t}$. The structure of DDPG is shown in Fig. \ref{f:DDPG}, where the local network will generate samples for the memory database. Then via randomly taking some samples from the memory database, the target critic and actor network will update ``softly" on the basis of value-iteration and policy-iteration, respectively. Note that, updating ``softly" means it will update these networks via discounting between the old and new target network. After that, it will transfer its network to the local network. 

According to the different functions of networks, we define them as
\begin{eqnarray}\label{eq:R2}
\text{local actor network:} & \mathbf{a} = \pi(\mathbf{s};\theta_l) \\\nonumber
\text{local critic network:} & \mathbf{R_t} = \mathbf{Q}(\mathbf{s};\mathbf{a};\omega_l) \\\nonumber
\text{target actor network:} & \mathbf{a} = \pi(\mathbf{s};\theta_t) \\\nonumber
\text{target critic network:} & \mathbf{R_t} =\mathbf{Q}(\mathbf{s};\mathbf{a};\omega_t) \nonumber
\end{eqnarray}
where $\pi$ represents the policy for selecting actions, $\mathbf{Q}$ represents the estimated Q value of different states and actions. In addition, $\theta_l, \theta_t, \omega_l, \omega_t$ represent the parameters of each network.
Then the performance objective function can be expressed as \cite{Philip2017Policy}
\begin{eqnarray}\label{eq:R3}
\mathbf{C}(\pi(\theta)) & =& \int_{\mathbf{s}}^{}\rho^{\pi}(\mathbf{s})\mathbf{Q}(s,\pi(\mathbf{s};\theta);\omega)d\mathbf{s} \\\nonumber
& =& \mathbf{E}_{\mathbf{s} \sim \rho^{\pi}}[\mathbf{Q}(\mathbf{s},\pi(\mathbf{s};\theta);\omega)]
\end{eqnarray}
where $\rho^{\pi}$ is the distribution function of $\mathbf{s}$, $\mathbf{C}(\pi(\theta))$ is the objective value of policy $\pi(\mathbf{s};\theta)$, and $\mathbf{E}_{\mathbf{s} \sim \rho^{\pi}}[\mathbf{Q}(\mathbf{s},\pi(\mathbf{s};\theta);\omega)]$ represents the mathematical expectation of $\mathbf{Q}$ value when the distribution of state $\mathbf{s}$ is $\rho^{\pi}$ and the corresponding actor policy is $\pi(\mathbf{s};\theta)$. Then according to the deterministic policy gradient theorem we can obtain \cite{David2014Deterministic}
\begin{eqnarray}\label{eq:R4}
\nabla_{\theta}\mathbf{C}(\pi(\theta)) \!\!&=\!\!& \int_{\mathbf{s}}^{}\rho^{\pi}(\mathbf{s})\nabla_{\theta}\pi(\mathbf{s};{\theta})\nabla_{\mathbf{a}}\mathbf{Q}(\mathbf{s},\mathbf{a};\omega)|_{\mathbf{a}=\pi(\mathbf{s};\theta)}d\mathbf{s} \nonumber\\
\!\!&=\!\!& \mathbf{E}_{\mathbf{s} \sim \rho^{\pi}}[\nabla_{\theta}\pi(\mathbf{s};\theta)\nabla_{\mathbf{a}}\mathbf{Q}(\mathbf{s},\mathbf{a};\omega)|_{\mathbf{a}=\pi(\mathbf{s};\theta)}] 
\end{eqnarray}

The updating of critic network and actor network based on \eqref{eq:R4} can be expressed as
\begin{eqnarray}\label{eq:R5}
&\delta_t = \mathbf{R}_t + \gamma \mathbf{Q}(\mathbf{s}_{t+1},\pi(\mathbf{s}_{t+1};\theta);\omega) - \mathbf{Q}(\mathbf{s}_t,\mathbf{a}_t;\omega) \\\nonumber
&\Delta\omega = \alpha_{\omega}\delta_t\nabla_{\omega}\mathbf{Q}(\mathbf{s}_t,\mathbf{a}_t;\omega) \\\nonumber
&\Delta\theta = \alpha_{\theta}\nabla_{\theta}\pi(\mathbf{s}_t;\theta)\nabla_{\mathbf{a}}\mathbf{Q}(\mathbf{s}_t,\mathbf{a}_t;\omega)|_{\mathbf{a}=\pi(\mathbf{s};\theta)} \nonumber
\end{eqnarray}
where $\delta_t$ is the temporal difference error, and $\gamma$ is the discount factor for calculating the objective function, $\alpha_{\omega}$ and $\alpha_{\theta}$ are the step-size for calculating the changing rate of $\omega$ and $\theta$, and the changing rate of $\omega$ and $\theta$ is $\Delta\omega$ and $\Delta\theta$. As for the HJB problem, only one step is considered, so $\delta_t = \mathbf{R}_t - \mathbf{Q}(\mathbf{s}_t,\mathbf{a}_t;\omega)$.

Finally, the updating of target networks based on \eqref{eq:R5} is expressed as
\begin{eqnarray}
\theta_t^n = \eta\theta_t + (1- \eta)\theta_t^n \nonumber\\
\omega_t^n = \eta\omega_t + (1- \eta)\omega_t^n \nonumber
\end{eqnarray}
where $\eta$ is the factor for updating the target networks, $\theta_t^n$ and $\omega_t^n$ represent the updated target networks. Then the DDPG can update iteratively. In the following subsection, more details of implementing DDPG in HJB are described below.

\subsubsection{Implementation of DDPG in HJB}
To implement the DDPG in HJB, we have designed new sampling rules, developed proper reward function, and determined suitable network architectures. The first step is to sample data from the memory database.
For the original DDPG, randomly selecting samples from the large memory database cannot work well for the HJB problem. Because the sensitivity of the initial adjoint variables, only a very small number of samples achieve good results. That means when integrating HJB equations from the randomly selected initial adjoint variables, only a few of them will reach the specified final boundary point, which makes DDPG diverge from the optimal solutions. 

To handle this issue, new database $\mathbf{D_g}$ is defined as: if and only if the actions $\mathbf{a}_t$ at state $\mathbf{s}_t$ can achieve reward $\mathbf{R}_t$ and satisfy $\mathbf{R}_t < \mathbf{R}_l$, then $\mathbf{a}_t \in \mathbf{D_g}$. Here, $\mathbf{R}_l$ is the upper bound of the reward. Then in each episode, a fixed part of samples will be generated from $\mathbf{D_g}$, and others will be randomly choosing from the  memory database. Additionally, $\mathbf{R}_l$ is set as a relative large number at the beginning, and gradually decreases along the iterations of DDPG. 
Therefore, $\mathbf{R}_l$ will be updated via $\mathbf{R}_l = \beta_l \mathbf{R}_l,\, \beta_l \in (0,1)$ when the average reward is less than $\mathbf{R}_l$.

Another issue of implementing DDPG is how to develop proper reward functions. When integrating HJB equations from the initial point, improper initial adjoint variables will make the integral of \eqref{eq:H2} significantly diverge from the specified final boundary value. To avoid this issue, the reward of $\lambda_0$ leading to an invalid final boundary value (value exceeds the allowable limit) is set to be $\mathbf{R}_t=N_1$. The reward of $\lambda_0$ leading to a large gap, e.g., $|\mathbf{x}(f)-\mathbf{x}_f|>L_1$ with $L_1$ being the upper bound of the gap value, but within the allowable limit is set to be $\mathbf{R}_t=N_2 + J$, where $J$ is defined in \eqref{eq:OCP1}. Additionally, when the integrated final boundary value from a given $\lambda_0$ and the specified one are close enough, e.g., $|\mathbf{x}(f)-\mathbf{x}_f| \leq L_2$ with $L_2$ set as a given threshold, the reward function is set to be $\mathbf{R}_t= J$. For the remaining cases, the reward function is set to be $\mathbf{R}_t=\beta_f|\mathbf{x}(f)-\mathbf{x}_f|+ J,\, \beta_f>0$. In summary, the reward function are determined by
\begin{eqnarray}\label{eq:R6}
\mathbf{R}_t=\begin{cases}
N_1, & \mathbf{x}(f) \text{ exceeds allowable limit}\\
N_2 + J, & |\mathbf{x}(f)-\mathbf{x}_f|>L_1\\
\beta_f|\mathbf{x}(f)-\mathbf{x}_f| + J, & L2 < |\mathbf{x}(f)-\mathbf{x}_f| \leq L_1\\
J, & |\mathbf{x}(f)-\mathbf{x}_f| \leq L_2
\end{cases}\nonumber
\end{eqnarray}

After that, a suitable network architectures is selected for DDPG. Via trial and errors, the final network architectures are shown in Table \ref{table_1}, where three layer neural networks are selected for both actor network and critic network. The activation function of their first layer and second layer are sigmoid function and hyperbolic tangent function (tanh), respectively. The output layer of actor network has tanh activation to generate continuous space and scaled by the same range of action space. Since the output of critic network is just the reward value, only one output is needed and linear activation is applied here. 
\begin{table}[!h]
	\caption{Architectures of actor and critic networks}
	\label{table_1}
	\begin{center}
		\begin{tabular}{|c||cc||cc|}
			\hline
			Iteration & ~~Actor&Network~~&~~Critic&Network~~\\
			\hline
			Layer & units& activation&units&activation\\
			\hline
			hidden 1 & 20*$dim_s$&sigmoid&20*$dim_s$&sigmoid\\
			\hline
			hidden 2 & 10*$dim_s$&tanh& 10*$dim_s$ & tanh\\
			\hline
			output & $dim_a$ &tanh& 1 & linear\\
			\hline
			
		\end{tabular}
	\end{center}
\end{table}

\section{SIMULATION RESULTS}
To verify the effectiveness of the proposed learning based HJB methods, a classical optimal control problem, the Brachistochrone problem, and a hypersensitive example are presented in this section. The results obtained from the proposed methods for both problems are analyzed and compared to the optimal solutions.

\subsection{The Brachistochrone Problem}
Given a fixed starting point $(x_0,y_0)$ and a terminal point $(x_f,y_f)$, a small bead can move in the gravity field with initial velocity $V_0 = 0$. The Brachistochrone problem is aiming to find the path that can achieve shortest time between the starting point and terminal point, which is formulated as
\begin{eqnarray}\label{eq:SR1}
&\min\limits_{\mathbf{u}} & \mathbf{J}_1 = \int_{0}^{t_f} 1dt \\
&\text{s.t.}  & \dot{\mathbf{x}} = V \cos u \nonumber\\
&             & \dot{\mathbf{y}} = V \sin u  \nonumber\\
&            & \mathbf{x}(0)=x_0,\, \mathbf{y}(0)=y_0,\, \mathbf{x}(f)=x_f,\, \mathbf{y}(f)=y_f \nonumber
\end{eqnarray}
where $V$ is the velocity of the small bead along the path, $u$ is the angle between the direction of $V$ and $x$-axis. As no external energy is added on the small bead, its dynamic energy can only be converted from its potential energy. Therefore, according to the energy conservation, $V$ can be expressed as
\begin{eqnarray}\label{eq:SR2}
V = \sqrt{2g\mathbf{y}}
\end{eqnarray}
By taking $\mathbf{x}$ as an independent variable, the problem is converted into
\begin{eqnarray}\label{eq:SR3}
&\min\limits_{\mathbf{u}} & \mathbf{J}_1 = \int_{x_0}^{x_f} \frac{1}{\sqrt{2g\mathbf{y}\cos u}}dx \nonumber\\
&\text{s.t.}  & \dot{\mathbf{y}} =\tan u \nonumber\\
&              &           \mathbf{y}(x_0)=y_0, \mathbf{y}(x_f)=y_f \nonumber
\end{eqnarray}
Then the HJB equations can be written as
\begin{eqnarray}\label{eq:SR3}
& \dot{\mathbf{y}} =\tan u \nonumber\\
& \dot{\lambda} = \frac{g}{V^3 \cos u} \nonumber\\
& \mathbf{y}(x_0)=y_0, \mathbf{y}(x_f)=y_f \nonumber
\end{eqnarray}

\begin{figure}[h]
	\centering
	\begin{subfigure}[h]{0.5\textwidth}
		\centering
		\includegraphics[width=8cm]{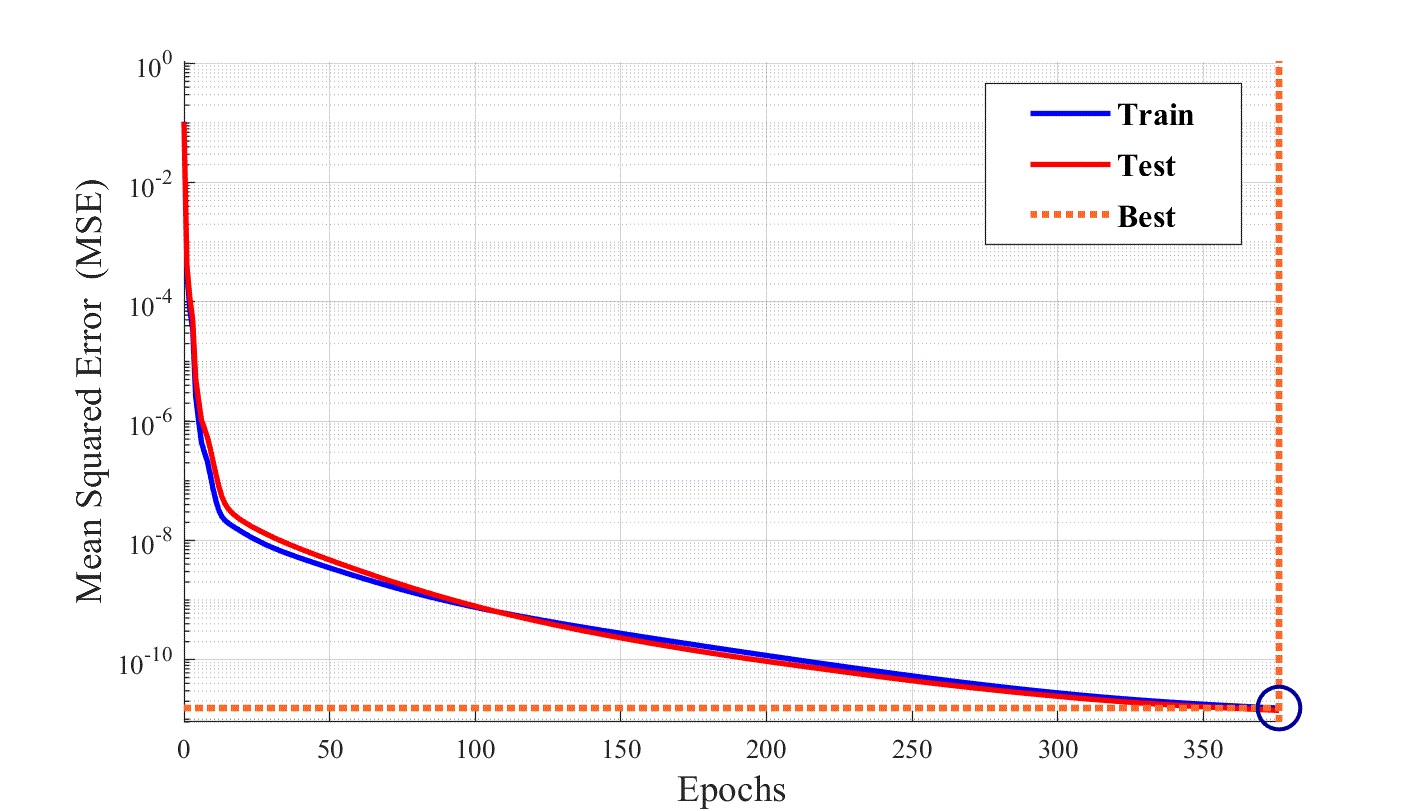}
		\caption{MSE between the output of neural network and tested data}
		\label{f:MSE}
	\end{subfigure}
	\begin{subfigure}[h]{0.5\textwidth}
		\centering
		\includegraphics[width=8cm]{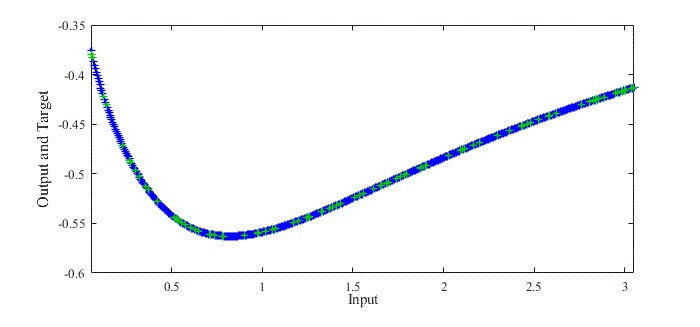}
		\caption{Comparison of predicted output from neural network and the optimal point, the blue line represents the predicted output of neural network and the green points represent the optimal point}
		\label{f:time}
	\end{subfigure}
	\begin{subfigure}[h]{0.5\textwidth}
		\centering
		\includegraphics[width=8cm]{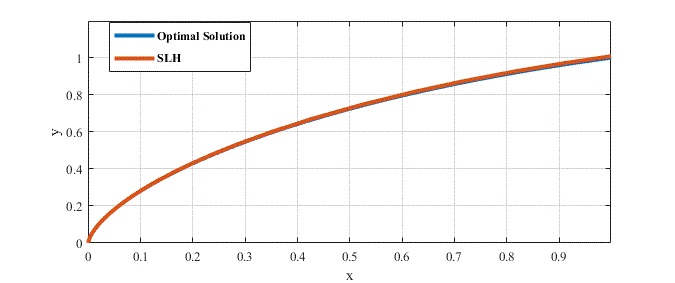}
		\caption{Comparison of the result of SLH with the optimal solution}
		\label{f:objective}
	\end{subfigure}
	\caption{The results of SLH for the Branchistochrone problem}
	\label{f:SLHA1}
\end{figure}

As the analytical solution of this problem is already achieved, the training database for SLH is obtained from its analytical solution. For the Brachistochrone problem, only two layers is included in the neural network of SLH, and the neurons of first layer is set as 10. Additionally, the activation function of first layer is the hyperbolic tangent sigmoid function (tansig). Via training the database, the learning result from SLH is shown in Fig. \ref{f:SLHA1}. The mean squared error (MSE) between trained neural network and the testing data is shown in Fig. \ref{f:MSE}. At 376 epochs, the MSE is $1.5398e^{-11}$ which means the neural network has converged and well trained. In Fig. \ref{f:time}, it can also be found that the outputs of neural network fit well with the testing data. In Fig. \ref{f:objective}, the result from SLH is compared to the optimal solution for one case starting from $(0,0)$ and ending at $(1,1)$. Result of this case shows that $ |\mathbf{x}(f)-\mathbf{x}_f| < 1e^{-3}$, which verifies that SHL achieves high precision when determining the initial adjoint variables. In conclusion, for the cases where validated HJB solutions are available to construct the training database, SLH can be a promising real-time method for solving these optimal control problems. In general, these validated HJB solutions can be obtained off-line from numerical methods, while this example uses analytical solutions to construct the training database.

\begin{figure}[!h]
	\centering
	\begin{subfigure}[h]{0.5\textwidth}
		\centering
		\includegraphics[width=8cm]{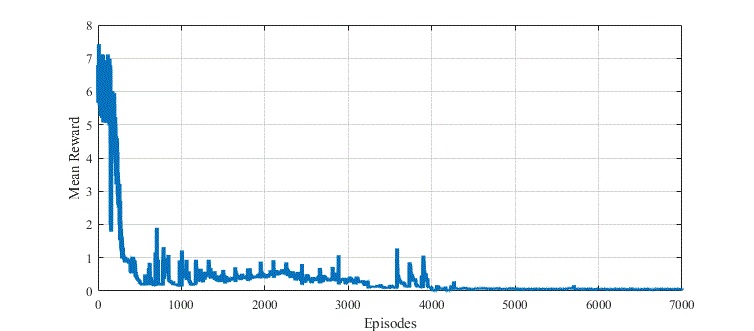}
		\caption{Mean reward of RLH for 7000 episodes}
		\label{f:mean}
	\end{subfigure}
	\begin{subfigure}[h]{0.5\textwidth}
		\centering
		\includegraphics[width=8cm]{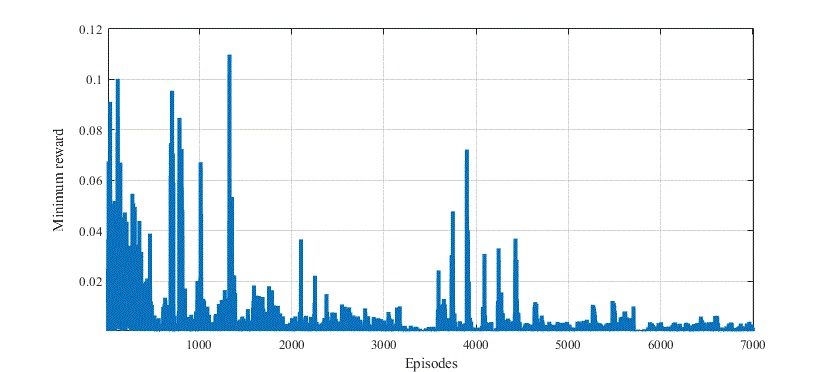}
		\caption{Minimum reward of RLH for 7000 episodes}
		\label{f:minimum}
	\end{subfigure}
	\begin{subfigure}[h]{0.5\textwidth}
		\centering
		\includegraphics[width=8cm]{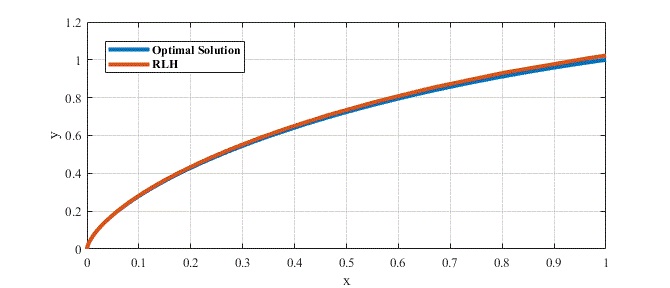}
		\caption{Comparison of the result of RLH with the optimal solution}
		\label{f:RLHA_c}
	\end{subfigure}
	\caption{The results of RLH for the Brachistochrone problem}
	\label{f:RLHA1}
	\vspace{-0.6cm}
\end{figure}

To apply RLH to solving the Brachistochrone problem, algorithm related parameters need to be properly selected. The maximum episodes is set to be 7000, the learning rate of actor network and critic network are set as $0.005$ and $0.01$, respectively. In addition, the total memory database size is set to be 30000 and the size of samples is 1000. The network architecture applied to this problem is shown in Table \ref{table_1}, and the soft updating factor for updating target network is $\eta = 0.01$. After off-line training, the results of RLH are shown in Fig. \ref{f:RLHA1}, which indicates that only 6000 episodes are required for RLH to converge. Additionally, the off-line training time is about 9200 seconds on a standard desktop. In Fig. \ref{f:RLHA_c}, the solution of RLH is compared to the optimal solution for the same case starting at $(0,0)$ and ending at $(1,1)$, which yields $ |\mathbf{x}(f)-\mathbf{x}_f| < 2e^{-2}$. These results verify that RLH can be effective for solving optimal control problems even in absence of reliable training dataset. Furthermore, comparing to the results of SLH, we find that SLH yields more precise results. This is due to the fact that SLH has labeled training database, which is more likely to lead to better results than just learning through randomly exploring samples. However, for cases when the proper labeled training database is not readily available, SHL is not applicable, while RLH becomes the only option. In conclusion, RLH can be a promising real-time method for solving optimal control even when the HJB solutions are difficult to be approximated via numerical methods. 

\subsection{The Hypersensitive Problem}
\begin{figure}[!h]
	\centering
	\begin{subfigure}[h]{0.5\textwidth}
		\centering
		\includegraphics[width=8cm]{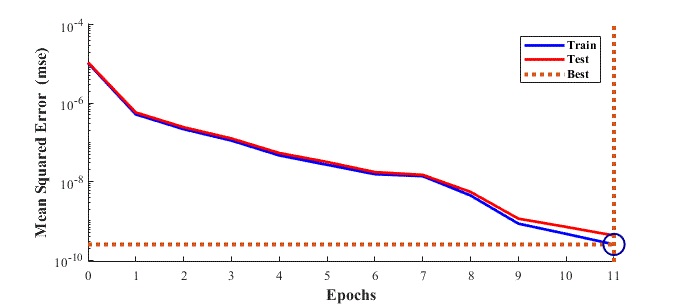}
		\caption{MSE between the output of neural network and tested data}
		\label{f:mean1}
	\end{subfigure}
	\begin{subfigure}[h]{0.5\textwidth}
		\centering
		\includegraphics[width=8cm]{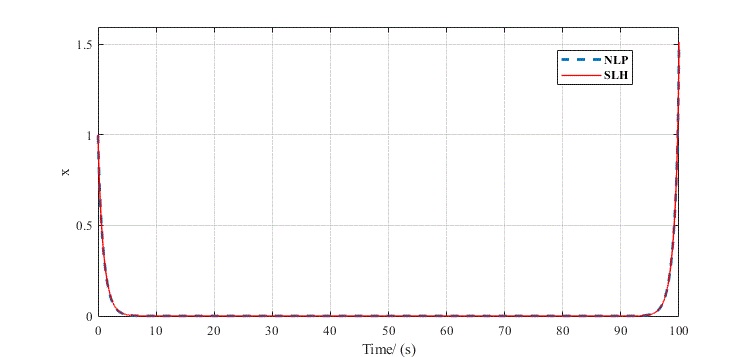}
		\caption{Comparison of the result of SLH with the solution from NLP}
		\label{f:minimum1}
	\end{subfigure}
	\caption{The results of SLH for the hypersensitive problem}
	\label{f:SLHA2}
\end{figure}

To verify the effectiveness of SLH for solving hypersensitive problems, an example is provided here,
\begin{eqnarray}\label{eq:M1}
&\min\limits_{\mathbf{u}} & \mathbf{J}_2 = \int_{0}^{t_f} (x^2+u^2)dt \\
&\text{s.t.}  & \dot{{x}} = -x^3 + u \nonumber\\
&            & \mathbf{x}(0)=1,\, \mathbf{x}(t_f)=1.5 \nonumber
\end{eqnarray}
Its HJB equations are expressed as
\begin{eqnarray}\label{eq:M2}
& \dot{\mathbf{x}} = -x^3 - \lambda/2, \\\nonumber
& \dot{\lambda} = -2x + 3x^2 \lambda\\\nonumber
& x(0)=1,\, x(t_f) = 1.5 \nonumber
\end{eqnarray}
To solve this problem, its stable and unstable phases are divided into 6 equal subsegments. For each subsegment, an independent artificial neural network is set to learn its initial adjoint variables. 
To verify the effectiveness of the proposed SLH for solving the hypersensitive problem, the result of SLH is compared with the solutions obtained from the nonlinear programming (NLP) solver, as shown in Fig. \ref{f:SLHA2}. Taking the training result of the first initial adjoint variables as example, the MSE between trained network and the testing data is shown in Fig. \ref{f:mean1}. At 11 epochs, the MSE is $2.8559e^{-10}$ which means the neural network has converged. By comparing the solution in Fig. \ref{f:minimum1}, it can be concluded that via learning multiple ``initial" adjoint variables for every subsegment, the integration to the final boundary point can be obtain with error $|\mathbf{x}(t_f)-1.5| < 1e^{-2}$. In addition, while solving the hypersensitive problem with NLP method, the longer the time interval is, the more discretized points are required, and the convergence of NLP method cannot be guaranteed. However, by off-line training, the real-time computational performance of SLH can always be guaranteed.

\section{CONCLUSIONS AND FUTURE WORKS}

In this paper, two learning based Hamilton-Jacobian-Bellman (HJB) approaches, named supervised learning Hamiltonian (SLH) and reinforcement learning Hamiltonian (RLH), are proposed to find the adjoint variables for HJB equations in real-time. By restoring validated solutions of HJB as training database, the supervised learning is able to learn the mapping function between the initial adjoint variables and boundary conditions off-line. Then the trained neural network can be applied to find the proper initial adjoint variables for the trained problem with new given boundary conditions in real-time. In addition, for hypersensitive HJBs, SLH can be extended by searching initial adjoint variables of multiple subsegments of the solution. In  absence of reliable training dataset where SLH is not applicable, the RLH is proposed to train without validated database. By setting proper neural network, reward function, and some super parameters, RLH is promising for solving optimal control problems in real-time.  


In the future, there are three directions to pursue. Firstly, as the hypersensitive problem is only solved via SLH in this paper, we will try to apply the RLH to solve hypersensitive cases as well. Secondly, as the accuracy of RLH is lower than that of SLH, we will try to improve the accuracy of RLH. Finally, many advanced technologies of reinforcement learning have been developed recently, we will integrate these technologies into RLH to further improve the computational performance.

%



\bibliographystyle{IEEEtran}
\bibliography{Refer}

\end{document}